\documentclass[11pt, a4paper, reqno]{amsart}

\usepackage{amsfonts,amsmath,amsthm}
\usepackage{amssymb}
\usepackage{latexsym}
\usepackage{calrsfs}

\usepackage[english]{babel}

\theoremstyle{plain}
\newtheorem{theorem}{Theorem}[section]
\newtheorem{theo}[theorem]{Theorem}

\newtheorem{cor}[theorem]{Corollary}

\theoremstyle{definition}

\theoremstyle{remark}

\def\stackrelunder#1#2{\mathrel{\mathop{#2}\limits_{#1}}}
\newcommand{\simsubset}{\stackrelunder{^\sim}{\subset}}

\newcommand{\eps}{\varepsilon}
\newcommand{\N}{{\mathbb N}}

\newcommand{\R}{{\mathbb R}}
\newcommand{\D}{{\mathbb D}}

\newcommand{\Id}{\mathrm{Id}}
\newcommand{\C}{\mathbb{C}}

\newcommand{\dopu}{{:}\allowbreak\ }

\newcommand{\rest}[2]{#1\raisebox{-0.3ex}{\mbox{$\mid_{#2}$}}}

\makeatletter
\def\@kcite#1#2{{%
  \@citestyle[\citeform{#1}\if@tempswa, K#2\fi]}}
\@ifundefined{kcite }{%
  \expandafter\let\csname kcite \endcsname\kcite
  \edef\kcite{\@nx\protect\@xp\@nx\csname kcite \endcsname}%
}{}
\def\kcite#1{{\@citestyle[K#1]}}%
\def\kcite{\cite}

\makeatother

\newcounter{abc}   
\newcounter{iiiii} 

\newenvironment{aequivalenz}
{\setcounter{iiiii}{0}
\begin{list}%
{{\rm (\roman{iiiii})}}
{\usecounter{iiiii}
\parsep=0pt plus 1pt
\topsep=1pt plus 2pt minus 1pt
\itemsep=1pt plus 2pt minus 1pt
\leftmargin=3\baselineskip \labelsep=.6\baselineskip
\labelwidth=2.4\baselineskip
\rightmargin 0pt}%
}
{\end{list}}

\newenvironment{statements}%
{\setcounter{abc}{0}
\begin{list}%
{{\rm (\alph{abc})}}
{\usecounter{abc}
\parsep=0pt plus 1pt
\topsep=1pt plus 2pt minus 1pt
\itemsep=1pt plus 2pt minus 1pt
\leftmargin=3\baselineskip \labelsep=.6\baselineskip
\labelwidth=2.4\baselineskip
\rightmargin 0pt}%
}
{\end{list}}

{\begin{list}%
{$\bullet$}
{\leftmargin=3\baselineskip
\labelsep=\baselineskip \labelwidth=2.5\baselineskip
\parsep=0pt plus 1pt
\topsep=1pt plus 2pt minus 1pt
\itemsep=1pt plus 2pt minus 1pt
\rightmargin 0pt}%
}
{\end{list}}

\makeatletter
\newcommand{\Version}{%
\renewcommand{\@oddfoot}{\small\hfil\texttt{Version of \today}}%
\renewcommand{\@evenfoot}{\small\texttt{Version of \today}\hfil{}}%
}
\newcommand{\Versionempty}{
\renewcommand{\ps@empty}{%
\renewcommand{\@oddhead}{}%
\renewcommand{\@evenhead}{}%
\renewcommand{\@evenfoot}{\small\hfil\texttt{Version vom \today}}%
\renewcommand{\@oddfoot}{\small\texttt{Version vom \today}\hfil{}}%
}%
}

\makeatother

\newcommand{\bea}{\begin{eqnarray*}}
\newcommand{\eea}{\end{eqnarray*}}
\newcommand{\beq}{\begin{equation}}
\newcommand{\eeq}{\end{equation}}
\newcommand{\begsta}{\begin{statements}}
\def\endsta{\end{statements}}
\newcommand{\begaeq}{\begin{aequivalenz}}
\def\endaeq{\end{aequivalenz}}

\newcommand{\pel}{Pe{\l}\-czy\'{n}\-ski}

\numberwithin{equation}{section}

\makeatletter
\relax 
\select@language{english}
\@writefile{toc}{\select@language{english}}
\@writefile{lof}{\select@language{english}}
\@writefile{lot}{\select@language{english}}
\@writefile{toc}{\contentsline {section}{\tocsection {}{1}{Introduction}}{1}}
\@writefile{toc}{\contentsline {section}{\tocsection {}{2}{$M$-ideals}}{2}}
\@writefile{toc}{\contentsline {section}{\tocsection {}{3}{Almost isometric embeddings}}{5}}
\@writefile{toc}{\contentsline {section}{\tocsection {}{4}{Unconditionality}}{8}}
\@writefile{toc}{\contentsline {section}{\tocsection {}{5}{Operators on $L_1$}}{11}}
\@writefile{toc}{\contentsline {section}{\tocsection {}{6}{Extensions of operators into $C(K)$-spaces}}{14}}
\@writefile{toc}{\contentsline {section}{\tocsection {}{}{References}}{16}}
\@writefile{toc}{\contentsline {section}{\tocsection {}{}{Papers by Nigel Kalton}}{16}}
\@writefile{toc}{\contentsline {section}{\tocsection {}{}{Other papers}}{17}}
\makeatother

\begin{document}

\title{Nigel Kalton's work in isometrical Banach space theory}

\author{Dirk Werner}

\date{March 14, 2011}



\thanks{This is the first draft of a paper written for the Nigel Kalton
  Memorial Website \texttt{http://mathematics.missouri.edu/kalton/},
  which is scheduled to go online in summer 2011.}

\address{Department of Mathematics, Freie Universit\"at Berlin,
Arnimallee~6, \qquad {}\linebreak D-14\,195~Berlin, Germany}
\email{werner@math.fu-berlin.de}


\maketitle

\thispagestyle{empty}

\section{Introduction}

Nigel Kalton was one of the greatest mathematicians of the last
40~years, although he did his best to conceal this fact. An outsider
wouldn't recognise the mathematical giant that he was in this modest
person who was always friendly and good-humoured and who was more than
willing to share his ideas with everyone. 

Nigel published more than 260 papers (including several books) not
only in Banach and quasi-Banach space theory, but also in so diverse
fields such as game theory, continued fractions, harmonic analysis,
operator semigroups  and convex geometry. 
Every single of these papers contains a deep contribution by Nigel,
often taking care of the most difficult case that his coauthors would
have to leave open without his help. He had a wide interest in
mathematics, and his problem solving abilities were legendary. For
instance, once after a colloquium talk on continued fractions he got
hooked on the subject and redeveloped the theory for himself over one
weekend, eventually solving the problem exposed in the talk. 

I met Nigel for the first time at the conference on Banach spaces in
Mons in 1987. It so happened that on the day after the conference we
were both waiting for the same train to Paris, but not for the same
coach: he told me that he always rides the first class, adding, ``I'm
snobbish.'' Of course he couldn't be more wrong in his
self-assessment! Some years later he solved a big problem in $M$-ideal
theory (see Section~\ref{sec2}), a problem, where we, a group of fresh
Ph.D.s in Berlin, couldn't get anywhere. He emailed me a file with his
solution, and this was the beginning of our collaboration, in which
more often than not I felt like a pedestrian next to a racing-car. 
In June 2010, after a talk at the conference in Valencia with a
somewhat set-theoretical flavour, I reminded him of the quote from
Star Trek, ``It's mathematics, but not as we know it.'' I knew that
this would strike his sense of humour; I did not know that this would
be the last time I saw him. 

In the next few sections I will try to survey some of Nigel's
contributions to Banach space theory. I will restrict myself to
problems of an isometric nature, but even this narrower area is still
so rich that omissions and misconceptions will be
inevitable. Certainly, the only way to do justice to Nigel's genius
would be to not only paraphrase the main results, but to expound all
the ideas contained in his papers. I have to leave this to an abler
mathematician. 

In the following, papers by Nigel will be cited in the form \kcite{Kal}
and other papers in the form \cite{LinPel2}. The bibliography  will
first list Nigel's 
cited papers chronologically, and then the other papers
alphabetically.

\section{$M$-ideals} \label{sec2}

An \textit{$M$-ideal} $V$ of a Banach space $E$ is a closed subspace such that
the dual admits an $\ell_1$-direct decomposition $E^*= V^\bot \oplus_1
W$ for some closed subspace $W\subset E^*$. (In other words, $V^\bot$
is an $L$-summand of $E^*$.) The notion of an $M$-ideal was
introduced by E.M.~Alfsen and E.G.~Effros \cite{AlEf}; for a detailed study
one may consult \cite{HWW}. Of course, $c_0$ is an $M$-ideal in
$\ell_\infty$, and Dixmier proved back in 1951 that $K(H)$, the space
of compact operators on a Hilbert space $H$, is an $M$-ideal in
$L(H)$, the space of bounded operators. By the end of the 1980s more
examples of Banach spaces $X$ for which $K(X)$ is an $M$-ideal in
$L(X)$ were known; basically, these examples were $\ell_p$-sums
(${p>1}$) or $c_0$-sums of finite-dimensional spaces and certain of
their subspaces and quotients. On the other hand, several necessary
conditions were known: $X$ has to be an $M$-embedded space (i.e., $X$
is an $M$-ideal when canonically embedded into $X^{**}$), which
implies for example that $X^*$ has the RNP, and $X$ must have 
the metric compact approximation property; this was proved by
P.~Harmand and \AA.~Lima \cite{HaLi}. 
For subspaces of $X\subset \ell_p$ the converse
was proved by C.-M.~Cho and W.B.~Johnson \cite{ChJo1}: the metric compact
approximation property is sufficient for $K(X)$ to be an $M$-ideal in
$L(X)$.   (Later Nigel obtained this in complete generality, see
Corollary~\ref{cor2.1.a}.)  

Let us recall what this approximation property means. 
A Banach space $X$ has the \textit{metric compact approximation property} if
there is a net of compact operators $K_i\dopu X\to X$ of norm ${\le1}$
that converges pointwise to the identity: $K_ix\to x$ for all $x\in
X$. Actually, an even stronger approximation property holds if $K(X)$
is an $M$-ideal in $L(X)$: one can achieve that $\limsup \|\Id
-2K_i\|\le1$ and both $K_i\to \Id$ and $K_i^*\to \Id$ pointwise
(\textit{shrinking unconditional metric compact approximation
  property}). Here the \textit{shrinking} bit  derives from the fact
that also $K_i^*x^*\to x^*$ for all $x^*$, like for the projections
associated with a shrinking basis, and the \textit{unconditionality}
is hidden in 
the norm condition  $\limsup \|\Id -2K_i\|\le1$; see the beginning of
Section~\ref{sec4}. 

Although the $M$-ideal problem for $K(X)$ was intensively studied in
the 1980s, there were no conditions on $X$ known that were necessary and
sufficient for $K(X)$ to be an $M$-ideal in $L(X)$; the best result by
then was given by W.~Werner \cite{WW1}: 
$K(X)$ is an $M$-ideal in $L(X)$ if and
only if $X$ has the metric compact approximation property by means of
a net $(K_i)$ satisfying 
$$
\limsup \|SK_i + T(\Id-K_i)\|\le \max\{\|S\|,\|T\|\} \qquad\forall
S,T\in L(X).
$$
This was a big achievement, but still the condition is so complicated
that one cannot check easily that it is fulfilled for $X=\ell_2$. 

This was the moment when Nigel got interested in the problem, probably
after some eventually successful bugging by Gilles Godefroy. (It
should be added that in the 1980s people working on $M$-ideals of
compact operators often used techniques from Nigel's much quoted paper
\kcite{Kal}.\footnote{It got quoted 47~times according to the Mathematical 
Reviews database;
but Nigel once told me that it was accepted for publication only at
the third attempt.})
In his paper \kcite{Kal-M} he offered an entirely new approach based on what he
called property~$(M)$ and property~$(M^*)$. Here are the definitions. 
A Banach space $X$ has \textit{property~$(M)$} if whenever $(x_i)$ is a
bounded weakly null net and $x,y\in X$ satisfy $\|x\|=\|y\|$, then 
$$
\limsup \|x_i+x\| = \limsup \|x_i+y\|, \eqno(1)
$$
and  $X$ has \textit{property~$(M^*)$} if whenever $(x_i^*)$ is a
bounded weak$^*$ null net and $x^*,y^*\in X^*$ satisfy $\|x^*\|=\|y^*\|$, then 
$$
\limsup \|x_i^*+x^*\| = \limsup \|x_i^*+y^*\|. \eqno(2)
$$
These notions can be recast in the language of types on Banach spaces;
see Section~\ref{sec6} below. 

It is easy to see that $(M^*)$ implies $(M)$; for the converse see
the remarks following Theorem~\ref{theo2.2}. 

Nigel's theorem is as follows. 

\begin{theo}\label{theo2.1}
The following assertions about a Banach space $X$  are equivalent:
\begaeq
\item
$K(X)$ is an $M$-ideal in $L(X)$.
\item
 $X$ has
property~$(M)$, does not contain a copy of $\ell_1$ and has the
unconditional metric compact approximation  property, i.e., there is a
net of compact operators satisfying $K_ix\to x$ for all $x$ and
$\limsup\|\Id-2K_i\|\le1$. 
\item
$X$ has property~$(M^*)$  and has the
unconditional metric compact approximation  property.
\endaeq
\end{theo}

Actually, in his paper Nigel only deals with the case of separable
spaces and the sequential versions of $(M)$ and $(M^*)$, but the
extension to the general case doesn't offer any difficulties; it can
be found in \cite[page~299]{HWW} for example. One should remark that
the sequential 
property~$(M)$ is trivially satisfied in Schur spaces (where by
definition weakly convergent sequences are norm convergent), but these
are excluded by the requirement that $\ell_1$ does not embed into $X$
in Theorem~\ref{theo2.1}.

Let me point out that it is trivial to verify that the $\ell_p$-spaces
for $1<p<\infty$ and in particular Hilbert spaces satisfy Nigel's
conditions; to see that $\ell_p$ has $(M)$ just note that for $\ell_p$
we have
$$
\limsup \|x+x_i\|^p = \|x\|^p + \limsup\|x_i\|^p 
\eqno(3)
$$
so that (1) becomes obvious. 

Using property~$(M)$ Nigel has been able to solve a number of open
problems in the theory of $M$-ideals. For example, he proved a general
version of the theorem of Cho and Johnson mentioned above: 

\begin{cor}\label{cor2.1.a}
If $K(X)$ is an $M$-ideal in $L(X)$ and $E\subset X$ has the
  compact metric approximation property, then $K(E)$ is an $M$-ideal
  in $L(E)$ as well.
\end{cor}

He also showed that Orlicz sequence spaces and more generally modular
sequence spaces can be renormed to have property~$(M)$; if in addition
such a space $X$ has a separable dual, then $X$ can be renormed so
that $K(X)$ becomes an $M$-ideal in $L(X)$. In the other direction, a
separable non-atomic order continuous Banach lattice can be renormed
to have property~$(M)$ if and only if it is lattice isomorphic
to~$L_2$. 

Another interesting corollary is that spaces $X$ with property~$(M)$ contain
subspaces isomorphic to $\ell_p$; more precisely, there exists $1\le
p<\infty$ so that $\ell_p$ embeds almost isometrically into $X$ (cf.\
Section~\ref{sec3} for this concept) or $c_0$ embeds almost
isometrically into $X$. The proof relies on a deep theorem due to
J.-L.~Krivine \cite{Kri}; a particular consequence is the theorem, originally
obtained by J.~Lindenstrauss and L.~Tzafriri, that an
ifinite-dimensional subspace of an Orlicz sequence space $h_M$
contains a copy of some $\ell_p$ or of $c_0$.

Concerning $L_p=L_p[0,1]$ it was known by 1990 that $K(L_p)$ is not an
$M$-ideal in $L(L_p)$ for $p\neq2$, and conversely that an $M$-ideal
in $L(L_p)$ is necessarily a two-sided ideal for $1<p<\infty$; indeed
this is so since $L_p$ and its dual are uniformly convex \cite{ChJo2}. But in
our joint paper \kcite{KalW1} 
it was shown that, for $p\neq1,2,\infty$, there are
no nontrivial $M$-ideals in $L(L_p)$ whatsoever (nontrivial meaning
different from $\{0\}$ and the whole space), and if $1<p,q<\infty$
then $K(\ell_p(\ell_q^n))$ is the only nontrivial $M$-ideal in
$L(\ell_p(\ell_q^n))$. In both these results we had to assume complex
scalars; the proofs use arguments involving hermitian operators. 

In the paper \kcite{AndCazKal}, with coauthors G.~Androulakis and
C.D.~Cazacu, Nigel 
takes his construction of spaces with property~$(M)$ still further in
that he considers Fenchel-Orlicz spaces \cite{Turett}. The definition of
these spaces is similar to that of Orlicz sequence spaces, but they
are built on a Young function on $\R^n$ rather than $\R$ and consist of
vector-valued sequences. Now Nigel and his coauthors proved that Fenchel-Orlicz
spaces can be renormed to have propert~$(M)$ and that many interesting
Banach spaces have a represention as a Fenchel-Orlicz space. This is
in particular so for the ``twisted sums'' $Z_p$, $1<p<\infty$, from
\kcite{KaltonPeck}. These are ``extreme'' counterexamples to the
three-space problem for $\ell_p$: $Z_p$ is not isomorphic to $\ell_p$,
yet contains a subspace $Y_p$ isomorphic to $\ell_p$ such that
$Z_p/Y_p$ is isomorphic to $\ell_p$ as well. In the language of
homological algebra, $Z_p$ is a nontrivial twisted sum of $\ell_p$
with itself, i.e., there is a short exact sequence $0\to \ell_p \to
Z_p\to \ell_p \to 0$ that does not split. Nigel has contributed a lot
to twisted sums, but this is another story. 

In a subsequent publication \kcite{Kal-Wer}  Nigel pursued
an idea mentioned at the end of \kcite{Kal-M}, namely to decide whether the
unconditionality assumption, i.e., $\limsup \|\Id-2K_i\|\le1$, in
Theorem~\ref{theo2.1} is actually 
needed. It turns out that this is not so. 

\begin{theo}\label{theo2.2}
For a separable Banach space $X$, $K(X)$ is an $M$-ideal in $L(X)$ 
if and only if $X$ has
property~$(M)$, does not contain a copy of $\ell_1$ and has the
 metric compact approximation  property.
\end{theo}

The proof of the if-part consists of two steps: first to show that $X$
has property~$(M^*)$ as well, which is much more difficult than the
implication $(M^*)$ $\Rightarrow$ $(M)$, and then to construct, using
property~$(M^*)$, from a compact approximation of the identity
satisfying $\limsup\|K_n\|\le1$ another compact approximation of the identity
satisfying $\limsup\|\Id-2L_n\|\le1$. This is done by a skipped
blocking decomposition argument. Meanwhile other and simpler arguments
for the second step have been given by \AA.~Lima \cite{Lima1995},
E.~Oja \cite{Oja2000}  and
O.~Nygaard and M.~P\~oldvere \cite{NygPold}. 

This theorem was proved while I was a visitor at the University of
Missouri in 1993. Let me commit myself to some personal recollections
at this stage. The day I arrived, Nigel asked me what I was working
on. One of the questions had to do with Banach spaces $X$ for which
$K(X\oplus_p X)$ is an $M$-ideal in   $L(X\oplus_p X)$, like
$X=\ell_p$. The conjecture was that such a space should be similar to
$\ell_p$, more precisely such an $X$ should embed almost isometrically
into an $\ell_p$-sum of finite-dimesional spaces. In fact, in the
paper \cite{Dirk4} I had previously formulated the bold conjecture that all
spaces for which $K(X)$ is an $M$-ideal in $L(X)$ are stable in the
sense of Krivine and Maurey \cite{KriMau} and that
one should be able to deduce from that that $K(X\oplus_p X)$ is an
$M$-ideal in   $L(X\oplus_p X)$ for some~$p$. The first half was
disproved by Nigel in \kcite{Kal-M} whereas he did prove the second part to be
correct. Almost on the spot he suggested an idea how to tackle the
problem. Of course, it took me some time to digest it, and after a
week or so I understood what he had in mind. I then suggested to use
an ultrapower argument at some stage of the proof, upon which Nigel
said, ``Oh, I missed that point!'' -- only to come up with a much
better idea that eventually solved the problem. I also pointed out a
relation to work by Bill Johnson and Morry Zippin, and Nigel asked, ``Can
they do it without the approximation property?'', which was the case,
and he added, ``Then we can do without the approximation property
too!'' I'll describe the outcome in the next section. 

\section{Almost isometric embeddings}\label{sec3}

Let us start with some vocabulary. We say that a (separable) Banach
space $X$ has property~$(m_p)$ if, whenever $x_n\to 0$ weakly,
$$
\limsup \|x+x_n\|^p = \|x\|^p + \limsup\|x_n\|^p \qquad\forall x\in X
$$
if $p<\infty$, resp. 
$$
\limsup \|x+x_n\| = \max\{ \|x\|,  \limsup\|x_n\| \} \qquad\forall x\in X
$$
for $p=\infty$. It is clear that $\ell_p$ has $(m_p)$ for $p<\infty$
(cf.\ (3) above) and that $c_0$ has $(m_\infty)$. If $K(X\oplus_p X)$
is an $M$-ideal, then $X$ has $(m_p)$, as proved by Nigel \kcite{Kal-M}. 

The Johnson-Zippin space $C_p$ is an $\ell_p$-sum of a sequence of 
finite-dimen\-sional spaces $E_1,E_2,\dots$ that are dense in all
finite-dimensional spaces with respect to the Banach-Mazur distance. A
Banach space $X$ embeds almost isometrically into $Y$ if for each
$\eps>0$ there is a subspace $X_\eps\subset Y$ such that
$d(X,X_\eps)\le 1+\eps$, $d(X,X_\eps)$ denoting the Banach-Mazur
distance. We use a similar definition for $X$ to be almost isometric
to a quotient of $Y$. Note that any two versions of $C_p$ (built on
different $E_k$) embed almost isometrically  into each other.  

In \kcite{Kal-Wer}, the following result is proved. 

\begin{theo}\label{theo3.1}
Suppose $X$ is a separable Banach space not containing $\ell_1$. Let
$1<p<\infty$. Then $X$ has $(m_p)$ if and only if $X$ embeds almost
isometrically into $C_p$. Likewise, $X$ has $(m_\infty)$ if and only
if $X$ embeds almost isometrically into $c_0$. 
\end{theo}

The proof uses again a skipped blocking decomposition
technique. 

Because of the duality of the property $(m_p)$, one obtains the
following corollary.

\begin{cor}\label{cor3.2}
Let $1<p<\infty$. If $X\subset C_p$, then $X$ is almost isometric to a 
quotient of $C_p$, and if $X=C_p/Z$, then  $X$ is almost isometric to a 
subspace of $C_p$. 
\end{cor}

This is an almost isometric refinement of a result due to Johnson and
Zippin who proved the corresponding isomorphic result \cite{JohZip2}. On the
other hand, isomorphic versions of Theorem~\ref{theo3.1} using tree
conditions were later obtained by Nigel for $p=\infty$ \kcite{Kal-QJM01} 
\label{page5}and E.~Odell and
Th.~Schlumprecht  for $p<\infty $ \cite{OdeSch}. 

In the context of $L_p$-spaces more can be proved. First of all, if
$X\subset L_p=L_p[0,1]$ has property~$(M)$, then it has $(m_r)$ for
some $r$; if $1<p\le2$, then $p\le r\le 2$, and if $2<p<\infty$, then
$r=2$ or $r=p$. We now have \kcite{Kal-Wer}:

\begin{theo}\label{theo3.3}
Suppose $1<p<\infty$, $p\neq2$, and let $X\subset L_p$ be
infinite-dimensional. Then the following are equivalent:
\begaeq
\item
$B_X$, the unit ball of $X$, is compact in $L_1$ (i.e., with respect
to the topology inherited from $L_1$).
\item
$X$ has property $(m_p)$.
\item
$X$ embeds almost isometrically into $\ell_p$.
\endaeq
If $p>2$, then \mbox{\rm{(i)--(iii)}} are equivalent to each of the
following:
\begaeq
\item[\rm(iv)]
$X$ is isomorphic to a subspace of $\ell_p$.
\item[\rm(v)]
$X$ does not contain a copy of $\ell_2$.
\endaeq
\end{theo}

Our proof was (iii) $\Rightarrow$ (ii) $\Rightarrow$ (i) $\Rightarrow$
(iii)  and (iii) $\Rightarrow$ (iv) $\Rightarrow$ (v)
$\Rightarrow$ (i) for $p>2$, where (iii) $\Rightarrow$ (iv)
$\Rightarrow$ (v) is trivial, as is (iii) $\Rightarrow$ (ii). 
Originally the equivalence of (iv) and
(v) for $p>2$ is due to Bill Johnson and Ted Odell \cite{JohOde}, see
also \cite{John3}. 

Concrete examples of subspaces of $L_p$ with $(m_p)$ are the Bergman
spaces consisting of all analytic functions on the unit disc
$\D=\{z\in\C\dopu |z|<1\}$ for
which $\int_\D |f(x+iy)|^p\,dx\,dy<\infty$. Likewise, the ``little''
Bloch space has $(m_\infty)$. In \kcite{Kal-Coll}, Nigel proved that
also the ``little'' 
Lipschitz space has $(m_\infty)$, thus showing that it is an $M$-ideal
in its bidual, a problem left open in \cite{BerWer}. 

Results concerning isometric embeddings of subspaces of $L_p$ into
$\ell_p$ were proved by F.~Delbaen, H.~Jarchow and A.~\pel\
\cite{DelJarPel}; as is 
often the case in isometric considerations in $L_p$, one has to
distinguish whether or not $p$ is an even integer. 

In a conversation in 1998, Nigel once suggested to prove a result
similar to Theorem~\ref{theo3.1} for property~$(M)$. His idea was to
show that such spaces embed into Fenchel-Orlicz spaces. I am sure that
he knew an outline of the argument, but as far as I know the proof has
never been written down. 

There is also an isomorphic version of Theorem~\ref{theo3.1} for
$p=\infty$, devised by Nigel and his coauthors G.~Godefroy and
G.~Lancien \kcite{GodKalLan}. 
The main result of their paper is that the Banach space
$c_0$ is determined by its metric, that is:

\begin{theo}\label{theo3.4}
If a Banach space $X$ is Lipschitz isomorphic to $c_0$, that is, if
there is a bijective map $T\dopu X\to c_0$ with $T$ and $T^{-1}$
Lipschitz, then $X$ is linearly isomorphic to $c_0$, that is, $T$ can
be chosen linear.  
\end{theo}

To prove this, they first show that $X$ embeds  isomorphically into
$c_0$, for which an isomorphic version of property~$(m_\infty)$ and of
Theorem~\ref{theo3.1} is needed. To conclude the proof of
Theorem~\ref{theo3.4} one has to appeal to known properties of
subspaces of $c_0$. Theorem~\ref{theo3.4} is one of the most
remarkable achievements in the nonlinear theory of Banach spaces. 

In their paper \kcite{GodKalLi1}, Nigel together with G.~Godefroy and D.~Li
addressed the problem of extending results like Theorem~\ref{theo3.3}
to the case of subspaces of~$L_1$. 
There are two intrinsic difficulties that do not occur in the case
$p>1$. For one thing, the Haar basis is unconditional in $L_p$ for
$p>1$, but not in $L_1$, and this was an essential ingredient in the
proof of Theorem~\ref{theo3.3}. Also, the $L_1$-topology on the unit ball of
a subspace $X\subset L_p$ is certainly locally convex, whereas the
$L_r$-spaces for $r<1$ are not locally convex and hence the  
$L_r$-topology on the unit ball of
a subspace $X\subset L_1$ need not be locally convex. (Here we enter
the world of quasi-Banach spaces, one of Nigel's favourite
areas.)
So in order to be able to study when subspaces of $L_1$ embed into
$\ell_1$, one has to assume some unconditionality. 

A separable Banach space has the \textit{unconditional metric
  approximation property} (UMAP for short) 
\label{page6}if there is a sequence of
finite rank operators such that $F_nx\to x$ for all $x$ and
$\|\Id-2F_n\|\to1$; the latter obviously implies that $\|F_n\|\to 1$
as well. We have already encountered a variant in
Theorem~\ref{theo2.1}; the definition of UMAP is due to Nigel and Pete
Casazza \kcite{CasKal}. We shall have more to say on this in
Section~\ref{sec4}.  

In the next theorem, one of the main results from \kcite{GodKalLi1}, $\tau_m$
denotes the topology of convergence in measure, i.e., the topology of
the $F$-space~$L_0$. 

\begin{theo}\label{theo3.5}
Let $X$ be a subspace of $L_1$ with the approximation property. The
following statements are equivalent.
\begaeq
\item
$X$ has the UMAP, and $B_X$ is relatively compact for the
topology~$\tau_m$. 
\item
$B_X$ is $\tau_m$-compact and $\tau_m$-locally convex. 
\item
For any $\eps>0$, there exists a weak\/$^*$ closed subspace $X_\eps$
of $\ell_1$ with Banach-Mazur distance $d(X,X_\eps)\le1+\eps$.
\endaeq
\end{theo}

By a result of Rosenthal, $B_X$ is $\tau_m$-relatively compact if and
only if $X$ fulfills a strong quantitative version of the Schur
property, called the $1$-strong Schur property. 

The paper \kcite{GodKalLi1} also contains a very interesting
counterexample as to 
possible generalisations of the previous theorem. 

\begin{theo}\label{theo3.6}
 There exists a subspace $X$ of $L_1$ with the approximation property,
 whose unit ball is $\tau_m$-compact but not $\tau_m$-locally
 convex. In particular, $X$ fails the UMAP.
\end{theo}

The construction of this space would not have been possible without
Nigel's insight into the nature of $L_p$-spaces for $0\le p<1$---for
most of us a no-go area---in particular the strange world of needle
points and the failure of the Krein-Milman theorem there, cf.\ \cite{Rob76}.

Another type of embedding result is contained in Nigel's work with
Alex Koldobsky \kcite{KalKol}; it concerns subspaces of the quasi-Banach space
$L_p$ for $p<1$. It is known that a Banach space $X$ embeds into some
$L_p$, $p<1$, isomorphically if and only if  $X$ embeds into all
$L_r$, $0<r<1$, isomorphically. (The embedding into $L_r$ for $r\le p$
is clear since $L_p$ embeds into those $L_r$ isometrically; the issue
is the range of $r$ between $p$ and~$1$.) Nigel proved in \kcite{Kal85}
that embedding into $L_1$ is equivalent to the
embeddability of $\ell_1(X)$ into $L_p$. 
As for the corresponding isometric question, A.~Koldobsky \cite{Kol96}
showed that there is a finite-dimensional Banach space that embeds
isometrically into $L_{1/2}$, but not into $L_1$. Using stable random
variables, Nigel and Alex Koldobsky obtain a vast generalisation.

\begin{theo}\label{theo3.7}
Let $0<p<1$; then there is an infinite-dimensional Banach space $E_p$
that embeds isometrically into $L_p$, but not into any $L_r$ for
${p<r\le 1}$. 
\end{theo}

The space $E_p$ has a representation as $\ell_1 \oplus_{N_p} \R$ by
means of some very cleverly chosen absolute norm $N_p$ on $\R^2$, and
the embedding of $E_p$ into $L_p$ is realised by means of a sequence
of independent random variables having a $1$-stable (i.e., Cauchy)
distribution. In addition to this example $E_p$, a second example is
constructed that is isomorphic to a Hilbert space, viz.\ $F_p= \ell_2
\oplus_{\tilde N_p} \ell_2$. 

\section{Unconditionality}\label{sec4}

We have already mentioned the notion of unconditional metric
approximation property (UMAP) on page~\pageref{page6}.
The UMAP was introduced in the paper \kcite{CasKal} by Nigel and Pete Casazza.
Let us explain what ``unconditional'' refers to here. Following Nigel
and Pete one can obtain from an approximating sequence $(F_n)$ with
$\|\Id-2F_n\|\to1$, for a given $\eps>0$, another approximating sequence
$(F'_n)$ such that for $A_n=F_n'-F_{n-1}'$ and all~$N$
$$
\biggl\|\sum_{n=1}^N \eps_n A_n \biggr\| \le 1+\eps
$$
whenever $\eps_n=\pm1$; this should be compared to the estimate
$$
\biggl\|\sum_{n=1}^N \eps_n e_n^*(x)e_n \biggr\| \le (1+\eps)\|x\|
$$
for $(1+\eps)$-unconditional bases. Hence the epithet ``unconditional.''
Replacing finite rank operator by compact operators in the
approximating sequence, one arrives at the notion of unconditional metric
compact approximation property (UMCAP). Apart form studying these
properties, \kcite{CasKal} also contains the proof of the following stunning
result concerning the ordinary metric approximation property (MAP). 

\begin{theo}\label{theo4.1}
If a separable Banach space has the MAP, then it even has the
commuting MAP, meaning that there are commuting finite rank operators
with $\|T_n\|\le1$ and $T_nx\to x$ for all~$x$.
\end{theo}

For reflexive spaces, the same type of conclusion is proved for the
UMAP in \kcite{CasKal} as well, but in a paper by Nigel and
G.~Godefroy \kcite{GodKal97}, 
the result was proved in full generality.

\begin{theo}\label{theo4.2}
If a separable Banach space has the UMAP, then it even has the
commuting UMAP.
\end{theo}

The key to the proof is to look at some norm-$1$ 
approximating sequence $(F_n)$
and to study the limiting projection
$$
Px^{**} = w^*\mbox{-}\lim_{n\to\infty} F_n^{**}x^{**}
$$
in the bidual and to prove that its range is weak$^*$ closed. The
paper \kcite{GodKal97} also contains 
an embedding result that is similar in spirit to Theorem~\ref{theo3.1}
and Theorem~\ref{theo3.3}: 

\begin{theo}\label{theo4.3}
A Banach space $X$ has the UMAP if and only if, for each $\eps>0$, $X$
embeds isometrically as a $(1+\eps)$-complemented subspace into a space with
a $(1+\eps)$-unconditional Schauder basis. 
\end{theo}

Another outgrowth of \kcite{CasKal} is the concept of a $u$-ideal that is
studied in detail in the influential paper \kcite{GKS} by Nigel, G.~Godefroy
and P.~Saphar. 
Let $X\subset Y$ be Banach spaces. $X$ is said to be a
\textit{$u$-summand} in $Y$ if there is a projection $P$ from $Y$ onto
$X$ with $\|\Id-2P\|=1$; equivalently, one may decompose $Y=X\oplus
X_s$ in such a way that $\|x+x_s\|=\|x-x_s\|$ whenever $x\in X$,
$x_s\in X_s$. 
An easy example of a $u$-decomposition is the decomposition of $f\in
C[-1,1]$ into its even and odd part. 
An important special case is when $Y=X^{**}$; for
example, the bidual of $L_1$ admits such a decomposition, which is
even $\ell_1$-direct: $(L_1)^{**}=L_1 \oplus_1 (L_1)_s$, the so-called
Yosida-Hewitt decomposition. (In technical terms, $L_1$ is an example
of an $L$-embedded space; see Chapter~IV in \cite{HWW}.) 
Likewise, $X$ is a \textit{$u$-ideal} in $Y$ if $X^\bot$, its
annihilator, is a $u$-summand in~$Y^*$. By definition, every $M$-ideal is a
$u$-ideal, but also every order ideal in a Banach lattice is a
$u$-ideal. When working with complex scalars, it is more appropriate
to replace the condition $\|\Id-2P\|=1$ by $\|\Id-(1+\lambda)P\|=1$
for all scalars of modulus~$1$; correspondingly, one then speaks of
\textit{$h$-summands} and \textit{$h$-ideals}. 
In these definitions, $u$ stands for
``unconditional'' and $h$ for ``hermitian.''

The extensive paper \kcite{GKS} contains a wealth of information concerning
$u$-ideals and $h$-ideals. I will mention only a few aspects. First
of all, there are certain similarities to $M$-ideals, for example, the
$u$-projection is uniquely determined, and a $u$-ideal that does not
contain a copy of $c_0$ is a $u$-summand. On the other hand, if $X$ is
a $u$-ideal in its bidual, then the $u$-projection need not be the
canonical projection from the decomposition $X^{***}=X^\bot \oplus
X^*$, 
as for $M$-embedded spaces. (An example is $X=L^1$.) Let us say that
$X$ is a strict $u$-ideal (in its bidual) in this case. 

One of the
main results in \kcite{GKS} characterises strict $u$-ideals by means of a
quantitative version of Pe{\l}czy\'nski's property~$(u)$. To explain
  this, some notation is needed. Let $x^{**}\in X^{**}$ be such that
there is a sequence in $X$ converging to $x^{**}$ in the weak$^*$
topology of $X^{**}$. Define
$$
\kappa_u(x^{**}) =
\inf \biggl\{ \sup_n \Bigl\| \sum_{k=1}^n \eps_k x_k
  \Bigr\| \dopu \eps_k=\pm1,\ x_k\in X,\
  x^{**}=w^*\mbox{-}\sum_{k=1}^\infty x_k \biggr\}
$$
and let $\kappa_u(X)$ be the smallest constant such that
$\kappa_u(x^{**})\le C \|x^{**}\|$ for all $x^{**}$ in the sequential
closure of $X$ in $(X^{**},w^*)$. $\kappa_h$ is defined in the same
way, replacing $\eps_k=\pm1$ with $|\eps_k|=1$, $\eps_k\in\C$. The
Banach space $X$ has property $(u)$ if $\kappa_u(X)<\infty$. 

\begin{theo}\label{theo4.4}
Suppose $X$ does not contain a copy of $\ell_1$. Then $X$ is a strict
$u$-ideal in $X^{**}$ if and only if $\kappa_u(X)=1$. 
\end{theo}

As for spaces of operators, Nigel and his collaborators obtain the
following results.

\begin{theo}\label{theo4.5}
Let $X$ be a separable reflexive Banach space. Then $X$ has the UMCAP if and
only if $K(X)$ is a $u$-ideal in $L(X)$. 
\end{theo}

Generally, the results for $h$-ideals are more satisfactory and
precise, since one can apply the powerful machinery of hermitian
operators; for example (the notion of complex UMCAP should be
self-explanatory): 

\begin{theo}\label{theo4.6}
Let $X$ be a complex Banach space with a separable dual. Then $X$ has
the complex  UMCAP if and
only if $K(X)$ is an $h$-ideal in $L(X)$ and $X$ is an $h$-ideal in $X^{**}$. 
\end{theo}

\begin{theo}\label{theo4.7}
Let $X$ be a separable complex Banach space with the MCAP. 
Then $K(X)$ is an $h$-ideal in its bidual $K(X)^{**}$ if and only if
$X$ is an $M$-ideal in $X^{**}$ and has the complex UMCAP. 
\end{theo}

In their paper \kcite{CowKal} (with S.R.~Cowell) Nigel takes up the question of
embedding into a space with an unconditional basis, as in
Theorem~\ref{theo4.3}, but without the assumption of the approximation
property. Such a result was given by W.B.~Johnson and B.~Zheng
\cite{JohZhe}, but now the aim is to find an (almost) isometric 
version. The key to this are the asymptotic unconditionality properties $(au)$
and $(au^*)$. 
A separable Banach space $X$ has \textit{property $(au)$} if 
$$ 
\lim (\|x+x_d\|-\|x-x_d\|) =0 
$$ 
whenever $(x_d)$ is a bounded weakly null
net; in \kcite{Kal-QJM01} Nigel referred to this property as ``$X$ is
of unconditional type.'' 
If $X^*$ is separable, it is equivalent to use weakly null
sequences instead, and one obtains a notion that has been known by the
acronym WORTH in the literature. Likewise, $X$ has \textit{property
  $(au^*)$} (previously  ``$X$ is of shrinking unconditional type'') if
$$ 
\lim (\|x^*+x_n^*\|-\|x^*-x_n^*\|) =0 
$$ 
whenever $(x_n^*)$ is a
(necessarily bounded) weak$^*$ null sequence; due to the weak$^*$
metrisability of the unit ball there is no need to consider nets
here. 
These notions are reminiscent of the properties $(M)$,
$(M^*)$ and $(m_p)$, and
$(au^*)$ has also been considered by \AA.~Lima \cite{Lima1995}
under the name
$(wM^*)$. 
In general $(au^*)$ implies $(au)$, and 
the converse is true under additional hypotheses. 
The main result of \kcite{CowKal}  says the following.

\begin{theo}\label{theo4.8}
A separable Banach space  $X$ has
property $(au^*)$ if and only if $X$ embeds almost isometrically into
a space $Y$ with a shrinking $1$-uncondi\-tional basis; if $X$ is
reflexive, $Y$ can be taken to be reflexive as well.
\end{theo}

I will briefly mention other contributions by Nigel on the topic of
unconditional bases, Schauder decompositions and expansions. In a
series of papers with P.~Casazza or F.~Albiac and C.~Ler\'anoz, 
Nigel investigated
uniqueness of unconditional bases in Banach or quasi-Banach spaces
(\kcite{CasKal96}, \kcite{CasKal98}, \kcite{CasKal99},
\kcite{AlbKalLer03}, \kcite{AlbKalLer04}). 
For example, in \kcite{CasKal99} it is proved that although the
Tsirelson space $T$ admits a unique unconditional basis, this is not
so for $c_0(T)$. In the paper \kcite{DefKal} with A.~Defant the question of
existence of an unconditional basis in the space $P(^m E)$ of
$m$-homogeneous on a Banach space $E$ is considered. Sean Dineen had
conjectured that $P(^m E)$ has an unconditional basis if and only if
$E$ is finite-dimensional. This conjecture is vindicated in
\kcite{DefKal}; 
the proof uses local Banach space theory and some greedy basis theory. 

See also Theorem~\ref{theo5.8} below for unconditional expansions in $L_1$.

\section{Operators on $L_1$}\label{sec5}

Nigel looked at operators on function spaces like $L_1$ or indeed
$L_p$ for $p\le1$ in a vast number of papers; I will report on a small
sample of his results.

In \kcite{Kal-Ind78} he devised a very useful representation
theorem for operators on $L_p$, $0\le p\le1$, by means of random
measures. For $p=1$ the result is as follows.

\begin{theo}\label{theo5.1}
If $T\dopu L_1[0,1]\to L_1[0,1]$ is a bounded linear operator, then
there are measures $\mu_x$ on $[0,1]$, with $x\mapsto \mu_x\in M[0,1]$
weak\/$^*$ measurable, such that
$$
(Tf)(x)= \int_0^1 f(s)\,d\mu_x(s) \quad a.e.
$$
Moreover,
$$
\|T\|= \sup_{\lambda(B)>0} \frac1{\lambda(B)} \int_0^1 |\mu_x(B)|\,dx.
$$
\end{theo}  

Decomposing the measures $\mu_x$ into their atomic and continuous
parts $\mu_x^a$ and $\mu_x^c$ allows to define the corresponding
operators
$$
(T^af)(x)= \int_0^1 f(s)\,d\mu^a_x(s), \quad
(T^cf)(x)= \int_0^1 f(s)\,d\mu^c_x(s) .
$$
This permits Nigel to derive the following variant of a result due to
P.~Enflo and T.W.~Starbird \cite{EnfStar}. 

\begin{theo}\label{theo5.2}
If $T^a\neq0$, there is a Borel set $B$ of positive measure such that
$\rest{T}{L_1(B)}$ is an into-isomorphism whose range is
complemented. 
\end{theo}

Let me remark that it is still an open problem whether a complemented
infinite-dimensional
subspace of $L_1$ is isomorphic to $L_1$ or $\ell_1$. However, Enflo
and Starbird  have proved that $L_p$ is primary for
$p\ge1$. This means that whenever $L_p$ is isomorphic to a direct sum
$X\oplus Y$, then $X$ or $Y$ is isomorphic to~$L_p$. Using his
representation theorem, Nigel extends this result to $p<1$ and obtains
a new proof for $p=1$. Thus:

\begin{theo}\label{theo5.3}
$L_p$ is primary for $0<p<\infty$.
\end{theo} 

The representation theorem~\ref{theo5.1} is also used in the 
paper \kcite{KalBeata} of Nigel and Beata Randrianantoanina, where
they show that a 
surjective isometry on a real rearrangement invariant space $X$ on
$[0,1]$ different from $L_2$ has the form $(Tf)(s)= a(s)f(\sigma(s))$;
if $X$ is not isomorphic to $L_p$ for any $1\le p\le \infty$, then in
addition $|a|=1$ a.e.\ and $\sigma$ is measure preserving. 

In \kcite{GodKalLi2}, Nigel and his coauthors G.~Godefroy and D.~Li take up
the line of reasoning based on Theorem~\ref{theo5.1} to obtain results
of a more isometric flavour, for example a quantitative version of a result
due to D.~Alspach \cite{Als2}.

\begin{theo}\label{theo5.4}
There is a function $\varphi$ with $\lim_{\alpha\to 0^+}
\varphi(\alpha)=1$ such that if $T\dopu L_1\to L_1$ satisfies
$$
\|f\|\le \|Tf\| \le \alpha \|f\|\qquad \forall f\in L_1,
$$
then there is an isometry $J\dopu L_1\to L_1$ such that $\|T-J\|\le
\varphi(\alpha)$. 
\end{theo}

In other words, there is always an into isometry close to a given into
near-isometry. 

The paper also contains an example, of a similar nature as that of
Theorem~\ref{theo3.6}, to show that this result does not extend to
(near-) isometries from subspaces of $L_1$ to~$L_1$. 

Specifically, Nigel and his coauthors investigate certain subspaces of
$L_1$ that they call \textit{small subspaces} and operators that they call
\textit{strong Enflo operators}. The latter means that in the
representation of $T$ by the measures $\mu_x$ one has
$\mu_x(\{x\})\neq0$ on a set of positive measure. A subspace $X$ of
$L_1$ is called small if the mapping $f\mapsto \rest fA$ from $X$ to
$L_1(A)$ is not surjective for any $A\subset [0,1]$ of positive
measure. In other words, the inclusion $B_{L_1(A)} \subset kB_X$ is
false whenever $\lambda(A)>0$ and $k>0$, where we consider
$L_1(A)\subset L_1[0,1]$ in a natural way. The following short-hand
notation is now handy. Say
$$
M \simsubset N
$$
for two subset of $L_1$ if there is a positive constant $k$ such that
$M\subset kN$, i.e., $M$ is absorbed by~$N$.  
Hence $X$ is not small if $B_{L_1(A)} \simsubset B_X$.  

It is worthwhile to equip $L_1$ with the topology $\tau_m$ of
convergence in measure, defined above Theorem~\ref{theo3.5}. 
For a subspace $X\subset
L_1$, let $C_X$ be the closure of $B_X$ in $L_1$ with respect to
$\tau_m$.  The subspace $X$ is called \textit{nicely placed} if
$B_X=C_X$, that is, if its unit ball is closed for the topology
$\tau_m$. By a theorem due to A.V.~Buhvalov and G.J.~Lozanovski\u\i\ 
\cite[page~183]{HWW}, $X$ is nicely placed if and only if the Yosida-Hewitt
projection associated to the decomposition
$$
(L_1)^{**}=L_1 \oplus (L_1)_s
$$
maps $X^{\bot\bot}$ onto $X$. Hence, $X\subset L_1$ is nicely placed
if and only if it is $L$-embedded, i.e., $X^{**}=X\oplus_1
X_s$. Nicely placed subspaces were introduced in \cite{G-bien} and studied
intensively in many papers, see Chapter~IV in \cite{HWW} for a survey.

Returning to the paper \kcite{GodKalLi2}, let us note the following result.

\begin{theo}\label{theo5.5}
A subspace $X$ of $L_1$ is small if and only if no strong Enflo
operator satisfies $T(B_{L_1})\simsubset C_X$. Consequently, a nicely
placed subspace is small if and only if no operator from $L_1$ to $X$
is a strong Enflo operator. 
\end{theo}

Since strong Enflo operators have a nonzero atomic part, they fix a
copy of $L_1$ by Theorem~\ref{theo5.2}; therefore a nicely placed
subspace of $L_1$ not containing $L_1$ is small.

Another feature of operators that are not strong Enflo operators is
the equation
$$
\|\Id+T\|=1+\|T\|, 
$$
known as the \textit{Daugavet equation}. By the above, this is
fulfilled for all operators valued in a small nicely placed
subspace. The Daugavet equation is one of the technical ingredients in
the proof of the main result of \kcite{GodKalLi2}:

\begin{theo}\label{theo5.6}
Let $X$ and $Y$ be small subspaces of $L_1$ and suppose that there is
an isomorphism $S\dopu L_1/X \to L_1/Y$ with
$\max\{\|S\|,\|S^{-1}\|\}<1+\delta< 1.25$. Then there is an
invertible operator $U\dopu L_1\to L_1$ such that
$\|U\|\,\|U^{-1}\|\le (1+\delta)/(1-25\delta)$ and
$d_{\mathcal{H}}(U(B_X),B_Y)\le 71\delta/(1-25\delta)$, where
$d_{\mathcal{H}}$ denotes the Hausdorff distance. 
\end{theo}

If $X$ and $Y$ are additionally assumed to be nicely placed, the
conclusion can be strengthened: If $\|S\|\,\|S^{-1}\|=1+\alpha<2$,
then $U$ above can be chosen to map $X$ onto $Y$ and 
$\|U\|\,\|U^{-1}\|\le (1+\alpha)/(1-\alpha )$. A particular
consequence in this setting is: If $L_1/X$ and $L_1/Y$ are isometric,
then so are $X$ and $Y$. 

Leaving the field of small subspaces let us turn to \textit{rich
  subspaces}. These were introduced by A.~Plichko and M.~Popov
\cite{PliPop}; later 
the definition was slightly modified \cite{KadSW2} to accommodate the general
setting of Banach spaces with the Daugavet property introduced in
\cite{KadSSW}.   
A Banach space $X$ has the \textit{Daugavet property} if the Daugavet
equation 
$$
\|\Id+T\|=1+\|T\| 
$$
holds for all compact operators $T\dopu X\to X$; examples include
$L_1[0,1]$, $C[0,1]$, $L_\infty[0,1]$, the disc algebra, $L_1/H^1$ and
many other function spaces, but also more pathological spaces in the
spirit of Theorem~\ref{theo3.6} \cite{KW}. A subspace $Y$ of $X$ is
called rich if 
whenever $Y\subset Z\subset X$, then $Z$ has the Daugavet
property. (There are other equivalent reformulations of this, see e.g.\
the survey \cite{Wer-IMB}.) 

I did some work on these notions with Vova Kadets and
other coauthors. Upon hearing about some of our results, Nigel
immediately contributed the following theorem, published in our joint
paper \kcite{KKW1}, showing that rich subspaces are indeed the largest possible
proper subspaces of~$L_1$. Recall that $C_X$ is the closure of $B_X$
for the topology of convergence in measure.

\begin{theo}\label{theo5.7}
A subspace $X\subset L_1$ is rich if and only if, for each
$1$-codimensional subspace $H\subset X$, $\frac12 B_{L_1}\subset
C_H$. On the other hand, if $rB_{L_1}\subset C_X$ for some
$r>\frac12$, then $X=L_1$. 
\end{theo}

One of the consequences of Nigel's representation theorem
(Theorem~\ref{theo5.1}) is not only the primariness of $L_1$, but more
generally that whenever $L_1$ is isomorpic to an unconditional
Schauder decomposition $X_1 \oplus X_2 \oplus \cdots$, then one of the
$X_k$ is isomorphic to $L_1$. In other words, this result ponders on
possible or impossible representations of $\Id=\sum T_n$ as a
pointwise unconditionally convergent series. In our joint paper \kcite{KKW2}
(with V.~Kadets) we investigate this question further. This paper
introduces a class $\mathcal C$ of operators related to the narrow
operators of Plichko and Popov \cite{PliPop} and to the not sign preserving
operators of H.P.~Rosenthal \cite{Ros-emb}, but I will skip the exact
definition and will only mention that compact operators belong to this
class. Anyway, here is the result. 

\begin{theo}\label{theo5.8}
Let $X$ be a Banach space and $T,T_n\dopu L_1\to X$ be bounded
operators such that $Tf=\sum_n T_nf$ unconditionally for each $f\in
L_1$. If each $T_n$ is in $\mathcal C$, then $T$ is in $\mathcal C$.
\end{theo}

Nigel's student R.~Shvydkoy had obtained a similar result for the
narrow operators in the case $X=L_1$ in his Ph.D. thesis. 

Theorem~\ref{theo5.8} can be considered as a generalisation of
A.~Pe{\l}czy\'nski's classical theorem that $L_1$ does not embed into
a space with an unconditional basis. It also contains an unpublished
(in fact unwritten, but occasionally quoted) result due to
H.P.~Rosenthal \cite{Ros-stop} as a special case:  
$L_1$ does not even sign-embed into
a space with an unconditional basis.

\section{Extensions of operators into $C(K)$-spaces}\label{sec6}

It is a classical fact that whenever $E\subset X$ are Banach spaces
and $T_0\dopu E\to L_\infty[0,1]$ is a bounded linear operator, then
there exists an extension $T\dopu X\to L_\infty[0,1]$ of the same
norm; however, there need not exist a bounded linear extension
whatsoever if $L_\infty[0,1]$ is replaced by $C[0,1]$. By constrast,
Joram Lindenstrauss's memoir \cite{Lin-Mem} presents a detailed study of those
range spaces for which compact operators are extendible; it turns out,
that, for all pairs $E\subset X$, every compact operator $T_0\dopu
E\to F$ admits, for every $\eps>0$, an extension $T\dopu X\to F$ of
norm $\|T\|\le (1+\eps)\|T_0\|$ if and only if $F^*$ is isometric to a
space $L_1(\Omega,\Sigma,\mu)$. In particular this is true for
$F=C[0,1]$.

More recently, the problem of extension of  bounded operators into $C(K)$ was
reexamined by Bill~Johnson and Morry Zippin \cite{JohZip95} who obtained
\hbox{$(3+\eps)$}-extensions for the pair $(E,\ell_1)$ provided $E$ is
weak$^*$ (meaning $\sigma(\ell_1,c_0)$-) closed. 

In a series of papers that are enormous in wealth of ideas, technical
mastery and also in size (\kcite{Kal-QJM01}, \kcite{Kal-NY07},
\kcite{Kal-Isr07}, \kcite{Kal-Ill08}) Nigel has contributed to this circle
of ideas. I will now describe just a few of his findings. Let us say,
following Nigel, that the pair $(E,X)$, with $E$ a subspace of the
separable space~$X$,  has the
\textit{$\lambda$-$C$-extension property} if, given a bounded linear
operator $T_0\dopu
E\to C(K)$ into some separable
$C(K)$-space, there is  an extension $T\dopu
X\to C(K)$ of norm $\|T\|\le \lambda\|T_0\|$. (The same problem can be
studied for Lipschitz maps where the Lipschitz constant plays the role
of the norm. This setting is investigated in detail in
\kcite{Kal-NY07}, \kcite{Kal-Isr07}.)
In this language the Johnson-Zippin result says that $(E,\ell_1)$ has
the $(3+\eps)$-$C$-extension property for every $\eps>0$ if $E\subset
\ell_1$ is weak$^*$ closed. 

Now, Nigel improves the $(3+\eps)$-bound to $(1+\eps)$, and he also
obtains a converse: If $(E,\ell_1)$ has the $\lambda$-$C$-extension
property for some $\lambda>0$ and additionally 
$\ell_1/E$ has an unconditional FDD,
then there is an automorphism $S\dopu \ell_1\to\ell_1$ such that
$S(E)$ is weak$^*$ closed. To prove this, he makes a detailed study of
spaces embedding into $c_0$ and presents the tree characterisation of
such spaces mentioned on page~\pageref{page5}. He also proves that if 
 $(E,\ell_1)$ has the $(1+\eps)$-$C$-extension
property for all $\eps>0$ (the ``almost isometric $C$-extension
property''), then $E$ is weak$^*$ closed. 

In \kcite{Kal-Ill08} Nigel goes on to develop an approach to the $C$-extension
property that is rooted in the theory of types, introduced by
J.L.~Krivine and B.~Maurey \cite{KriMau}. Let $X$ be a separable Banach space. A
\textit{type} generated by a bounded sequence $(x_n)$ (or more
generally a bounded net) is a function of the form
$$
\sigma\dopu X\to \R,\quad \sigma(x)= \lim_{n\to\infty} \|x+x_n\|
$$
(provided all these limits exist); it is called a weakly null type if
$x_n\to0$ weakly. Property~$(M)$ can now be rephrased by saying that
for a weakly null type $\sigma$, $\sigma(x)$ depends only on $\|x\|$,
and $X$ has $(m_p)$ whenever each weakly null type has the form
$\sigma(x)= (\|x\|^p + \lim_n \|x_n\|^p)^{1/p}$. Property~$(au)$ is a
symmetry condition, viz.\ $\sigma(x)=\sigma(-x)$. Now $X$ is said to
have \textit{property~$(L)$} if two weakly null types coincide once they
coincide at~$0$:
$$
\sigma_1(0)= \sigma_2(0) \quad\Rightarrow\quad 
\sigma_1(x)= \sigma_2(x) \mbox{ for all }x.
$$
Property $(L^*)$ is defined similarly using weak$^*$ null types in the
dual. Clearly, $\ell_p$ has properties $(L)$ and $(L^*)$ for
$1<p<\infty$, and so do certain renormings of Orlicz sequence spaces
or Fenchel-Orlicz spaces. These renormings are similar in spirit to
the renormings that yield property~$(M)$,
but not identical. Indeed, if $X$ has both properties $(M)$
and $(L)$ and fails to contain $\ell_1$, then it has $(m_p)$ for some
$1<p\le\infty$. 

With the help of types, Nigel is able to characterise the almost
isometric $C$-extension property. I won't give the exact formulation
of this result, but will note only one special case.

\begin{theo}\label{theo6.1}
If $X$ is a separable Banach space with $(M^*)$ or $(L^*)$ and if
$E\subset X$, then $(E,X)$ has the almost
isometric $C$-extension property.
\end{theo}

In particular, this applies to (renormings of) the twisted sums
$Z_p$. Nigel also provides a particular renorming of $\ell_2$ to show
that the existence of almost isometric extensions need not guarantee
isometric extensions. 

Moreover, Nigel gives the first examples of spaces with the following
universal property: Whenever $E$ embeds into a separable space $X$
isometrically, then the pair $(E,X)$ has the  almost
isometric $C$-extension property.
Indeed, all weak$^*$ closed subspaces of $\ell_1$ have this property
as do the spaces of Theorem~\ref{theo3.6}. The proof again uses the
theory of types. 

An intriguing study of the corresponding isomorphic property is
contained in \kcite{Kal-Ill08}. Let us say
 that the pair $(E,X)$ of separable spaces has the
\textit{$C$-extension property} if, given a bounded linear
operator $T_0\dopu
E\to C(K)$ into some separable $C(K)$-space, there is  a bounded linear 
extension $T\dopu X\to C(K)$; here  it is enough to
consider $K=[0,1]$ by Milutin's theorem. A Banach space $E$ has the 
\textit{universal $C$-extension property} if $(E,X)$ has the
$C$-extension property whenever $E$ embeds into  $X$
isometrically. 

Nigel proves the following result:

\begin{theo}\label{theo6.2}
A separable Banach space $E$ has the universal $C$-extension property
if and only if $E$ is $C$-automorphic. 
\end{theo}

The latter means that whenever $E_1\subset C[0,1]$ and $E_2\subset
C[0,1]$ are isomorphic to $E$, then there is an automorphism $S\dopu
C[0,1]\to C[0,1]$ mapping $E_1$ to $E_2$. One can paraphrase this by
saying that there is essentially only one way to embed $E$ into
$C[0,1]$. 

It is a classical result due to J.~Lindenstrauss   and A.~\pel\ \cite{LinPel2}
that $c_0$ is $C$-automorphic. 
Theorem~\ref{theo6.2} enables Nigel to show that $c_0(X)$ is
$C$-automorphic as well if $X$ is (for example $X=\ell_1$), 
but $\ell_p$ is not $C$-automorphic for $1<p<\infty$.
 Indeed, for a certain superreflexive $Z\supset \ell_p$ with an
 unconditional basis, the pair $(\ell_p,Z)$ fails the $C$-extension
 property. If, however, $Z\supset \ell_p$ is a UMD-space with an
 unconditional basis, then $(\ell_p,Z)$ satisfies the $C$-extension
 property; in Nigel's own words, ``the appearance of the
 UMD-condition is quite mysterious.'' 

A technical device to prove
 these theorems are homogeneous mappings $\Phi\dopu X^* \to Z^*$ that are
 weak$^*$ continuous on bounded sets such that $\Phi(x^*)$ extends
 $x^*$ with a bound $\|\Phi(x^*)\|\le \lambda \|x^*\|$, for all
 $x^*\in X^*$. Such mappings were introduced by  M.~Zippin \cite{Zip03}
and are called Zippin selectors by Nigel.
 Based on this notion, the following
 technical key result on $X=\ell_p$ or, more generally, an
 $\ell_p$-sum of finite-dimensional spaces, $1<p<\infty$, is proved,
 where again types play an essential role:

\begin{theo}\label{theo6.3}
Let $X$ be as above.
 For a separable superspace $Z\supset X$, 
the pair $(X,Z)$ has the
 $C$-extension property if and only if $Z$ can be renormed so as to
 contain $X$ isometrically and such that 
$$ 
\lim_{n\to\infty}  \|z+x_n\| \ge 
\lim_{n\to\infty} (\|z\|^p+ \|x_n\|^p)^{1/p} 
$$ 
for all
 $z\in Z$ and all weakly null sequences $(x_n)\subset X$, provided
 both limits exist.  
\end{theo}

Many of Nigel's extension theorems have counterparts for Lipschitz
functions; for a detailed study see \kcite{Kal-NY07} and \kcite{Kal-Isr07}.



\def\refname{Papers by Nigel Kalton}
\makeatletter
\renewcommand{\@biblabel}[1]{[K#1]}
\makeatother

\def\refname{Other papers}
\makeatletter
\renewcommand{\@biblabel}[1]{~[#1]}
\makeatother

\end{document}